\documentclass[12pt]{article}

\usepackage{a4wide}   
\usepackage{amsmath}   
\usepackage{amssymb}   
\usepackage{stmaryrd}  
\usepackage{latexsym}  
\usepackage{amsxtra}   
\usepackage{amstext}   
\usepackage{bm}        
\usepackage{amsthm}    
\usepackage{amscd}     
\usepackage[mathscr]{eucal}

\usepackage{epsfig}    
\usepackage{bbm}       

\usepackage[latin1]{inputenc}

\usepackage{epic}      
\usepackage{eepic}     
\usepackage{calc}      

\usepackage{xy}        
\xyoption{all}

\cleardoublepage 

\numberwithin{equation}{section}

\newcounter{counter}
\theoremstyle{plain} \numberwithin{equation}{subsection}

\newtheorem{thm}[equation]{Theorem}
\newtheorem{lem}[equation]{Lemma}

\newtheorem{prop}[equation]{Proposition}

\theoremstyle{definition}

\newtheorem{expl}[equation]{Example}

\theoremstyle{remark}

\theoremstyle{plain} \theoremstyle{definition}

\newcommand{\cI}{\mathcal{I}}

\newcommand{\cO}{\mathcal{O}}
\newcommand{\bC}{\mathbb{C}}
\newcommand{\bF}{\mathbb{F}}
\newcommand{\bH}{\mathbb{H}}
\newcommand{\bN}{\mathbb{N}}
\newcommand{\bP}{\mathbb{P}}

\newcommand{\bR}{\mathbb{R}}
\newcommand{\bS}{\mathbb{S}}

\newcommand{\bZ}{\mathbb{Z}}

\newcommand{\gK}{\mathbf{K}}
\newcommand{\gM}{\mathbf{M}}
\newcommand{\gO}{\mathbf{O}}
\newcommand{\gS}{\mathbf{S}}
\newcommand{\gT}{\mathbf{T}}
\newcommand{\gU}{\mathbf{U}}

\newcommand{\actson}{\curvearrowright}

\newcommand{\colim}{\operatorname{colim}}
\newcommand{\Diff}{\operatorname{Diff}}

\newcommand{\gaminf}{\Gamma_{\infty}}

\newcommand{\id}{\operatorname{id}}
\newcommand{\im}{\operatorname{im}}

\newcommand{\K}{\gK}
\newcommand{\ord}{\operatorname{ord}}

\newcommand{\Sl}{\operatorname{SL}}

\newcommand{\Sp}{\operatorname{Sp}}
\newcommand{\SU}{\operatorname{SU}}
\newcommand{\Th}{\operatorname{Th}}
\newcommand{\Td}{\operatorname{Td}}
\newcommand{\td}{\operatorname{td}}

\newcommand{\Tors}{\operatorname{Tors}}
\newcommand{\bTh}{\mathbb{T}\mathbf{h}}

\newcommand{\MU}{\gM \gU}
\newcommand{\infloop}{\Omega^{\infty}}
\newcommand{\loopinf}{\Omega^{\infty}}
\newcommand{\susp}{\Sigma^{\infty}}
\newcommand{\cpinf}{\bC \bP^{\infty}}
\newcommand{\MT}{\gM \gT \gS \gO(2)}
\newcommand{\Q}{\infloop \susp}

\begin{document}

\title{The low-dimensional homotopy of the stable mapping class group}

\author{Johannes F. Ebert\footnote{Supported by a fellowship within
the Postdoc-Prgramme of the German Academic Exchange Service (DAAD)}\\
e-mail: ebert@maths.ox.ac.uk}

\maketitle

\begin{abstract}
Due to the deep work of Tillmann, Madsen, Weiss and Galatius, the cohomology of the stable mapping class group $\gaminf$ is known with rational or
finite field coefficients. Little is known about the integral cohomology. In this paper, we study the first four cohomology groups. Also, we
compute the first few steps of the Postnikov tower of $B \gaminf^+$, the Quillen plus construction applied to $B \gaminf$. Our method relies on
the Madsen-Weiss theorem, a few known computations of stable homotopy groups of spheres and projective spaces and on a certain action of the
binary icosahedral group on a surface. Using the latter, we can also describe an explicit geometric generator of the third homotopy group $\pi_3
(B \gaminf)$.
\end{abstract}

\section{Introduction and overview}

Let $\Gamma_{g,b}$ be the group of isotopy classes of orientation-preserving diffeomorphisms of a connected compact surface $F_{g,n}$ of genus
$g$ with $b$ boundary components (both, diffeomorphisms and isotopies are assumed to fix the boundary pointwise), in other words, the
\emph{mapping class group}. We write $\Gamma_g:= \Gamma_{g,0}$. There are stabilization maps defined for $b >0$:

\begin{itemize}
\item $h: \Gamma_{g,b} \to \Gamma_{g,b+1}$ (gluing in a pair of pants along one boundary component),
\item $i:\Gamma_{g,b+1} \to \Gamma_{g+1,b}$, (glueing in a pair of pants along two boundary components),
\item $j:\Gamma_{g,b} \to \Gamma_{g,b-1}$, (glueing in a disc).
\end{itemize}

The homological stability theorem of Harer \cite{har}, as improved by Ivanov \cite{Iv}, \cite{Iv1}) asserts that $H_{k} (h;\bZ)$ and $H_{k}
(i;\bZ)$ are isomorphisms whenever $g \geq 2k+1$, while $H_{k} (j;\bZ)$ is an isomorphism when $g \geq 2k+2$. It follows that $H_k (B
\Gamma_{g,n})$ does not depend on $g$ or $n$ \emph{in the stable range}, i.e. as long as $g \geq 2k+2$. Let $\Gamma_{\infty,b}:= \colim ( \ldots \to
\Gamma_{g,b+1} \to Gamma_{g+1, b+1} \to \ldots$, where the composition $ i \circ h$ is used to form the colimit. The homology groups of this infinite mapping class group are the stable homology groups of the mapping class group. Due to the work of
Tillmann \cite{Till}, Madsen, Weiss \cite{MW} and Galatius \cite{Gal}, the stable homology of $\Gamma_g$ is completely understood, at least in
theory. Let us describe this story briefly. Even before Harers work, Powell \cite{Pow} showed that $\Gamma_g$ is perfect if $g \geq 3$, i.e. $H_1(B
\Gamma_g; \bZ)=0$. The Quillen plus construction from algebraic K-theory provides (see \cite{Ros}) a simply-connected space $B
\Gamma_{\infty,b}^{+}$ and a homology equivalence $B \Gamma_{\infty,b} \to B \Gamma_{\infty,b}^{+}$. The glueing maps $j$ induce a homotopy
equivalence $B \Gamma_{\infty,b}^{+} \simeq B \Gamma_{\infty}^{+}$. While $B \gaminf$ is by definition aspherical, $B \gaminf^+$ has many
nontrivial homotopy groups. Tillmann proved in \cite{Till} that $\bZ \times B \Gamma_{\infty}^{+}$ has the homotopy type of an infinite loop space. Later, she
and Madsen \cite{MT} constructed a spectrum $\MT$ and maps $\alpha_{g,b}: B \Gamma_{g,b} \to \infloop_0 \MT$. Madsen and Weiss showed that
$H_k(\alpha_{g,b})$ is an isomorphism if $g \geq 2k+2$, in other words, the limit map $B \Gamma_{\infty}^{+} \to \infloop_0 \MT$ is a homotopy
equivalence. One consequence of these results is the existence of a map
$\alpha_g: B \Gamma_g \to \infloop_0 \MT \simeq B\Gamma_{\infty}^{+}$ (without boundaries!) which induces isomorphisms in homology in the stable range.\\
The rational homology of $\Omega^{\infty}_{0} \MT$ is easy to determine (it is isomorphic to the rational homology of $BU$, see \cite{MT}). The
homology with $\bF_p$-coefficients for all primes $p$ was computed by Galatius \cite{Gal}. His
result is very complicated and it is very hard to destill explicit information form it, even in small dimensions.\\

The aim of this note is the study of the homotopy and homology groups of $B \gaminf^{+}$ up to dimension $4$. The first few homotopy groups of $\loopinf \MT$ and hence also of $B \gaminf^{+}$ were already computed in \cite{MT}. The result is 

\begin{thm}\label{homotopy}
$$ \pi_k(B \gaminf^+) \cong
\begin{cases}
0 & \text{if $k =1$}\\
\bZ & \text{if $k=2$}\\
\bZ/24 & \text{if $k =3$}\\
\bZ  & \text{if $k =4$.}\\
\end{cases}$$
\end{thm}

This follows from Serre's classical computation of the first values of $\pi_k (\Q\bS^0)$ and from a computation by Mukai \cite{Muk}.\\

\begin{thm}\label{homology}
The first few homology groups of the stable mapping class group are given by

$$ H_k(B \gaminf^+; \bZ) \cong
\begin{cases}
0 & \text{if $k =1$}\\
\bZ & \text{if $k=2$}\\
\bZ/12 & \text{if $k =3$}\\
\bZ^2  & \text{if $k =4$.}\\
\end{cases}$$
\end{thm}

The first homology was computed by Powell, as mentioned above, the second by Harer \cite{Har2}. Harer also proved that $H_3(B \gaminf)$ is a torsion group, see \cite{Har3} and determined the rank of $H_4(B \gaminf)$ (unpublished). Our method for the proof of theorem \ref{homology} is a quite indirect in the sense that we use Theorem \ref{homotopy} for it.
The relation between the homology and homotopy groups is expressed by the first stages of the Postnikov tower of $B \gaminf^{+}$.
The \emph{Postnikov tower} of a \emph{simple} space $X$ tells us how $X$ is built from the Eilenberg-MacLane spaces $K(\pi_n(X),n)$. We will compute the
Postnikov tower of $B \gaminf^+$ up to degree $4$. We introduce the following notation. A map $f: X \to Y$ is a \emph{$\pi_n$-isomorphism}, if
it induces an isomorphism $\pi_n (X) \to \pi_n(Y)$. A space $Y$ is \emph{$N$-coconnected} if $\pi_k(Y)=0$ for all $k > n$. If $X$ is a simple
space, then there is an $n$-coconnected space $\tau_{\leq n} X$ and a map $a_n : X \to \tau_{\leq n} X$, which is a $\pi_k$-isomorphism for $ k
\leq n$. These data are unique up to homotopy. The maps $a_n$ assemble to a homotopy-commutative diagram

$$\xymatrix{
                                                     & \ldots \ar[d]\\
                                                      &  \tau_{\leq n+1} X \ar[d]^{p_{n+1}}  \\
X \ar[dr]^{a_{n-1}} \ar[r]^{a_{n}} \ar[ur]^{a_{n+1}} \ar[uur] \ar[ddr] &  \tau_{\leq n} X \ar[d]^{p_{n}} \\
                                                      &  \tau_{\leq n-1} X \ar[d]\\
                                                      & \ldots
}
$$

We can assume that the maps $p_n$ are fibrations. The homotopy type of the space $X$ can be recovered as the inverse limit of this tower of fibrations. The homotopy
fiber of $p_n$ is an Eilenberg-MacLane space $K(\pi_n(X),n)$. Thus the fibration $p_n$ is classified by a map $k^n:  \tau_{\leq n-1} X \to
K(\pi_{n}(X);n+1)$ in the sense that $p_n$ is the homotopy-pullback of the path-loop fibration $K(\pi_n(X),n ) \to P K(\pi_n(X),n+1) \to
K(\pi_n(X),n)$ via $k^n$. Or, equivalently, $\tau_{\leq n}X$ is the homotopy fiber of the map $k^n$. Now we state the main result of this paper.

\begin{thm}\label{Postnikov}
The Postnikov invariants $k^2$ and $k^4$ of $B\gaminf^+$are trivial. The Postnikov-invariant $k^3:\tau_{\leq 2} B \gaminf^+ \simeq K(\bZ;2)
\simeq \cpinf \to K(\bZ/24;4)$ has order $2$ when considered as an element in $H^4(\cpinf; \bZ/24) \cong \bZ/24$.
\end{thm}

Given Theorem \ref{homotopy}, the theorems \ref{homology} and \ref{Postnikov} are equivalent - this follows easily from the Leray-Serre
spectral sequence. The first sentence of Theorem \ref{Postnikov} follows from the existence of a $\pi_2 $-isomorphism $\phi_2:B \gaminf^+ \to \cpinf$ and a
$\pi_4$-isomorphism $\phi_4:B \gaminf^+ \to K(\bZ;4)$. The existence of $\phi_2$ is easy because $\pi_1(B \gaminf^+)=0$, but we show the
following. The mapping class group $\Gamma_g$ acts on the cohomology of the surface $F_g$ preserving the intersection form, providing a group
homomorphism $\Gamma_g \to \Sp_{2g}(\bZ) \to \Sp_{2g}(\bR)$. Because $U(g) \subset \Sp_{2g}(\bR)$ is a maximal compact subgroup, $B \Sp_{2g}(\bR) \simeq B U(g)$. We obtain a map $\Phi_g: B \Gamma_g \to BU(g)$. This map classifies a complex vector bundle $V_1$,
which is also easy to describe. Consider the universal surface bundle\footnote{By definition, $\Gamma_{g}^{1}:=\pi_0 (\Diff(F_g,p)$, where $p \in F_g$. It is quite obvious that the forgetful map $B \Gamma_{g}^{1} \to B \Gamma_{g}$ is homotopy equivalent to the universal surface bundle on $B \Diff(F_g)$, as long as $g \geq 2$.} $\pi:B\Gamma_{g}^{1} \to B \Gamma_g$. One can choose a complex structure
on the vertical tangent bundle of this bundle, turning $\pi:B\Gamma_{g}^{1} \to B \Gamma_g$ into a family of Riemann surfaces. Let $V_1$ be the
vector bundle whose fiber over $x \in B \Gamma_g$ is the dual of the $g$-dimensional vector space of holomorphic $1$-forms on $\pi^{-1}(x)$. It
is a standard result that this is indeed a vector bundle once the topology is appropriately chosen \cite{AS}. It is easy to see that the
isomorphism class of this bundle does not depend on the choice of the complex structure and that it is classified by $\Phi_g$. It was shown in
\cite{MT} that there is a map $\Phi:B \gaminf^+ \to BU$ such that $\Phi_g = \Phi \circ \alpha_g$. The map $\Phi$ can be used to obtain
cohomology classes of $B \gaminf^+$. We set $\zeta_i := \Phi^* c_i$ and $\gamma_i := \Phi^* s_i$; $s_i \in H^{2i}(BU;\bZ)$ the universal integral Chern character class. It was shown by Morita \cite{Mor} that
$\gamma_{2i}$ is a torsion class.

\begin{prop}
The composition

$$\xymatrix{ B \gaminf^+ \ar[r]^{\Phi} & BU \ar[r]^-{c_1} & K(\bZ; 2)}$$

is a $\pi_2$-isomorphism. Consequently, $H^2(\gaminf^+; \bZ)\cong \bZ$ is generated by $\gamma_1$.
\end{prop}

This was already proven by Harer \cite{Har2} using different methods. To carry out the necessary computation, we use a
bordism-theoretic interpretation of $\pi_i(\infloop \MT)$ and the Hirzebruch-Riemann-Roch theorem. For the existence of the announced
$\pi_4$-isomorphism we use the Segal splitting of $\Q \cpinf$ \cite{Seg}.\\
To get the information about the third Postnikov invariant, we first construction an action of the binary icosahedral group $\hat{G}$, the fundamental
group of the famous Poincar\'e sphere, on the surface $F_{14}$. In other words, we construct a homomorphism $\rho:\hat{G} \to \Gamma_{14}$, which
yields a map $B \rho:B \hat{G} \to B \gaminf^+$. The next step is easily derived from the Lefschetz fixed point formula.

\begin{prop}\label{order24}
The characteristic class $B \rho^* \zeta_2 \in H^4(B \hat{G};\bZ) \cong \bZ/120$ has order $24$.
\end{prop}

This result, together with $\pi_3(B \gaminf^+) \cong \bZ/24$, suffices to finish the proof of theorem \ref{Postnikov} and thus also of
\ref{homology}.\\
Finally, we use the action $\hat{G} \actson F_{14}$ to construct a generator of $ \pi_3 (B \gaminf^+)$. This goes as follows. The group
$\hat{G}$ is perfect and acts freely on $\bS^3$. The quotient $M:=\bS^3/ \hat{G}$ is the Poincar\'e sphere. This gives a map $M \to B
\Gamma_{14}$. We can apply the plus construction to it and obtain a map $f: M^+ \to B \gaminf^+$. Because $M^+ \simeq \bS^3$, we have an element $\theta \in \pi_3 (B \gaminf^+)$.\\
The second main result of this note is

\begin{thm}\label{generator}
The element $\theta \in \pi_3(B \gaminf^+) \cong \bZ/24$ is a generator 
\end{thm}

The structure of the paper is as follows. In section \ref{even}, we recall the proof of Theorem \ref{homotopy} and prove the first sentence of Theorem \ref{Postnikov}. The section \ref{icosahedral} contains the construction of the $\hat{G}$-action on which the determination of the third Postnikov invariant, which is carried out in section \ref{end}, depends. The proof of Theorem \ref{generator} is also given in section \ref{end}.
\textbf{Acknowledgements:}
This paper is an improved version of a chapter of the author's PhD thesis \cite{Eb}. The author wants to express his thanks to his PhD advisor, Carl-Friedrich B\"odigheimer, for his constant support and patience. Also, I want to thank the Max-Planck-Institute for Mathematics in Bonn for financial support during my time as a PhD student. The final writing of this paper was done while the author stayed at the Mathematical Institute of the Oxford University, which was made possible by a grant from the Deutscher Akademischer Austauschdienst.\\
Finally, I want to thank Jeffrey Giansiracusa for teaching me Theorem \ref{homozopie}.

\section{The first few homotopy groups of $\loopinf \MT$}\label{even}

Given some classical computations of homotopy groups, it is not hard to compute the first five homotopy groups of $\loopinf \MT$. Recall that
there exists a fiber sequence of spectra \cite{GMT}

\begin{equation}\label{ravenel}
\xymatrix{ \MT \ar[r]^{\omega} & \susp \cpinf_{+} \ar[r]^{trf_{\bS^1}} & \Sigma^{\infty} \bS^{-1}. }
\end{equation}

The second map is the so-called \emph{circle transfer}. We are going to use the calculation of the circle transfer on low-dimensional homotopy
groups by Mukai \cite{Muk} and the classical low-dimensional computations of Serre. The first homotopy groups of the spectrum
$\Sigma^{\infty} \bS^0$ are (these values are tabularized in \cite{Tod})

\begin{tabular}[ht]{|c|c|c|c|c|c|c|}
\hline k & 1 & 2 & 3 & 4 & 5 & 6\\
\hline $\pi_{k}(\Q \bS^0)$ & $\bZ/2$ & $\bZ/2$ & $\bZ/24$ &$ 0$ & $0$ & $\bZ/2$. \\
\hline
\end{tabular}

The nontrivial element in $\pi_1$ is denoted $\eta$ (the Hopf element). We also need to know the effect of the multiplication by $\eta$ on
homotopy. $\eta^2$ is nonzero and $\eta^3$ has order $2$. The first values of $\pi_k
(\Q \cpinf)$ (without additional basepoint) are given by (\cite{Muk}, p.199)

\begin{tabular}[ht]{|c|c|c|c|c|c|c|c|}
\hline k & 1 & 2 & 3 & 4 & 5 & 6 & 7\\
\hline $\pi_{k}(\Q \cpinf)$ & $0$ & $\bZ$ & $0$ &$ \bZ$ & $\bZ/2$ & $\bZ$ & $\bZ/2$. \\
\hline
\end{tabular}

There is a product decomposition $\Q \cpinf_+ \simeq \Q \cpinf \times \Q \bS^0$. The effect of the circle transfer on second summand of $\pi_n(\Q \cpinf_+) = \pi_n(\Q \cpinf) \oplus \pi_n(\Q \bS^0)$ is given by the multiplication with $\eta$. Using this and Theorem 1.1 of \cite{Muk} we obtain that $\pi_k(\susp \cpinf_{+}) \to
\pi_{k+1}(\Q \bS^0)$ is surjective if $k \leq 5$. Thus the long exact homotopy sequence associated to \ref{ravenel} breaks into short pieces

$$0 \to \pi_k (\loopinf \MT) \to \pi_k (\Q \cpinf_{+}) \to \pi_{k+1} (\bS^0) \to 0 $$

for $k \leq 4$. An easy diagram chase finishes the proof of Theorem \ref{homotopy}.

\subsection{Geometric interpretation}

The spectra in the fiber sequence \ref{ravenel} are all Thom spectra and therefore the long exact homotopy sequence has a manifold-theoretic
interpretation\footnote{The author wants to thank Jeff Giansiracusa for explaining this interpretation to him.}. The interpretation of the
homotopy groups of a general Thom spectrum is a standard result. Let $X_{*}:=(X_0 \subset X_1 \subset X_2 \subset \ldots)$ be a sequence of
topological spaces and let $V_* \to X_*$ be a stable vector bundle of dimension $d \in \bZ$, i.e. a sequence $V_n \to X_n$ of real vector bundles of
dimension $n+d$ ($d$ can be negative; then $X_n$ is forced to be empty if $d+n <0$), together with isomorphisms $\epsilon_n: V_n \oplus \bR
\cong V_{n+1}|_{X_{n-1}}$. An important example of a stable vector bundle is the stable normal bundle $\nu_M$ of a closed smooth $m$-manifold
$M$; it has dimension $-m$. The \emph{Thom spectrum} $\bTh (V_*)$ it defined as follow: the $n$th space is $\bTh (V_*)_n:= \Th(V_n)$, the Thom
space of $V_n$ (the Thom space of the empty vector bundle over the empty space is the one-point space). The structural map $e_n : \Sigma
\Th(V_n) \to \Th(V_{n+1})$ is defined by the vector bundle map $\epsilon_n$.\\
The Pontryagin-Thom construction \cite{Rud} gives the interpretation of the homotopy groups $\pi_n (\bTh(V_*))$ (which agree with $
\pi_n(\loopinf \bTh(V_*))$ for $n \geq 0$). The group $\pi_n (\bTh(V_*))$ is isomorphic to the group of bordism classes of pairs $(M,\phi)$,
where $M$ is a closed smooth $n-d$-dimensional manifold and $\Phi: \nu_M \to V_*$ is a map of stable vector bundles. Because the suspension
spectrum $\susp X_+$ is the Thom spectrum of the trivial stable vector bundle, whose $n$th term is $X \times \bR^n$, it follows that $\pi_n
(\susp \cpinf_+)$ is the group of bordism class of triples $(M,b, E)$, $M$ a closed smooth $n$-manifold, $b$ a framing of $M$ and $E \to M$ a
complex line bundle. Consider the following stable vector bundle $-L$ of dimension $-2$. Its $2n$th term is the orthogonal complement
$L_{n}^{\bot}$ of the tautological line bundle $L_n \subset \bC \bP^{n-1} \times \bC^n$. Its $2n+1$st term is $L_{n}^{\bot} \oplus \bR$. The
Madsen-Tillmann spectrum $\MT$ is just the Thom spectrum of this stable vector bundle.\\
It follows that $\pi_n(\MT)$ is the group of bordism classes of triples $(M, W, f)$, $M$ a closed smooth $n+2$-manifold, $W \to M$ a complex
line bundle and $f: \bR^n \oplus W \cong TM$ is isomorphism of stable vector bundles. Equivalently, we can say that we have triples $(M,W,
g)$ with a stable vector bundle isomorphism $g: \nu_M \oplus W \cong \bR^{-n}$ instead. Given the description of the homotopy groups as bordism
groups and the maps in the fiber sequence \ref{ravenel}, it is quite easy to derive the geometric meaning of the maps in the long exact homotopy
sequence.

\begin{thm}\label{homozopie} (J. Giansiracusa, unpublished)
The maps in the long exact homotopy sequence of \ref{ravenel} are the following.
\begin{itemize}
\item The connecting homomorphism $\pi_{n+1}(\bS^{-1}) \to \pi_n (\MT)$ maps a framed manifold $(M^{n+2},b)$ to $(M, \bC, b \oplus
\id_{\bC})$ (second description).
\item The homomorphism $\pi_n (\MT) \to \pi_n (\susp \cpinf_{+})$ is represented
by the following procedure. Let $(M^{n+2}, L, f)$ as above be given. Choose a section $s$ of $L$ which is transverse to the zero section. Then
$s^{-1}(0)$ is an $n$-manifold and $f$ defines a framing of it.
\item The circle transfer $\pi_n (\susp \cpinf_+)$ maps a triple $(M,b,E)$ to the total space $\bS(E)$ of the sphere bundle of $E$, endowed with
the induced framing.
\end{itemize}
\end{thm}

We use Theorem \ref{homozopie} to give a description of the second and of the fourth homotopy group of $\MT$. There is a map of spectra $\lambda:\MT \to \Sigma^{-2}\MU$. It arises from the interpretation of $-L$ as a stable complex vector bundle of real dimension $-2$. The
effect of $\lambda$ on homotopy groups is the following. If $(M; L; f)$ represents an element in $\pi_n (\MT)$, then $f$ induces a stable bundle
isomorphism $\nu_M \cong \bR^n - L$, which gives a stable complex structure. We can thus view $M$ as a stable complex manifold and hence as an
element in $\pi_{n+2}(\MU)$. Recall also the Conner-Floyd map $\mu:\MU \to \K$ \cite{CF}. Bott periodicity gives an
isomorphism $\beta:\pi_{2n}(\K) \to \bZ$. Then $\beta \circ \mu$ sends $ [M] \in \pi_{2n}(\MU)$ to the Todd genus $\Td(M):= \langle \td (TM), [M]
\rangle \in \bZ$. The first components of the Todd class $\td(V)$ of a complex vector bundle are \cite{Hirz}: $\td_1(V)=\frac{1}{2} c_1(V)$,
$\td_2(V)=\frac{1}{12}(c_2(V) + c_{1}^{2}(V))$, $\td_3(V)=\frac{1}{24} (c_2(V) c_1(V))$.

\begin{thm}
The composition $\pi_2 (\MT) \to \pi_4 ( \MU) \to \pi_4 ( \K) \cong \bZ$ is an isomorphism.
\end{thm}

\textbf{Proof:} We already mentioned that $\pi_2 (\MT)  \cong \bZ$, but we need a generator. Look at the sequence

$$\xymatrix{
0 \ar[r] & \pi_2 (\MT) \ar[r] \ar@{=}[d]& \pi_2 (\susp \cpinf_{+}) \ar[r] \ar@{=}[d]& \pi_2 (\bS^{-1}) \ar[r]
\ar@{=}[d]& 0 \\
0 \ar[r] & \bZ \ar[r] & \bZ \oplus \bZ/2 \ar[r] & \bZ/24 \ar[r] & 0.}
$$

Generators of $\pi_2 (\susp \cpinf_{+})$ are the class of $\bS^2$ with the Hopf bundle and the canonical framing (order $\infty$) and the torus
with the Lie framing and the trivial line bundle (order $2$). The image of a generator of $\pi_2 (\MT)$ in $\pi_2 (\susp \cpinf_{+})$ is twelve
times the first generator plus the second one. It follows that a generator of $\pi_2 (\MT)$ can be described by the following data. It also
follows that these data exist.

\begin{itemize}
\item $M$ is a closed smooth $4$-manifold.
\item $L \to M$ is a complex line bundle.
\item There is a stable vector bundle isomorphism $\bC  \oplus L \cong TM$.
\item There is a section $s$ in $L$, transverse to the zero section, such that the zero set $S := s^{-1}(0)$ is an oriented surface with
self-intersection number $12$.
\end{itemize}

We claim that the Todd genus of such a manifold is necessarily equal to $1$. The Todd class of a $4$-manifold like $M$ is $\frac{1}{12}(c_2 (TM)
+ c_{1}(M)^{2}$. The third condition above ensures that $c_2 (TM) =0$ and hence $\td(TM)= \frac{1}{12}  c_1 (L)^2$ Thus $\Td (M)= \frac{1}{12}
\langle c_1 (L)^2; [M] \rangle = \frac{1}{12} 12$ by the fourth condition. \qed

A similar argument shows

\begin{prop}
The composition $\bZ \cong \pi_4 (\MT) \to \pi_6 ( \MU) \to \pi_6 ( \K) \cong \bZ$ is zero.
\end{prop}

We can also find a map $\pi_4$-isomorphism $\Omega^{\infty} \MT \to K(\bZ;4)$. This needs some preparation. The tautological line
bundle $L \to \cpinf$ defines a spectrum cohomology class in $\K^0 (\cpinf_{+})$ and thus a spectrum map $\Sigma^{\infty} \bC
\bP^{\infty}_{+} \to \K$. It induces a map $L: \Q \cpinf_{+} \to \Omega^{\infty} \K = \bZ \times BU$. It is a theorem by Segal \cite{Seg}
that there exists a splitting $\Phi: \bZ \times BU \to \Q \cpinf_{+}$, such that $L \circ \Phi \sim \id$ (it is not an infinite loop
map). Moreover, the fiber of $L$ has finite homotopy groups. Because $\pi_4 (\Q \cpinf_{+}) \cong \bZ$ and because $\pi_4 (BU) \cong
\bZ$, it follows that $L$ is a $\pi_4$-isomorphism. Moreover an inspection of \ref{ravenel} shows that $\omega:\MT \to \Q \cpinf_{+}$
is a $\pi_4$-isomorphism as well.

\begin{prop}
Let $a,b \in \bZ$. The effect of the cohomology class $ac_{1}^{2} + b c_2$, considered as a map $BU \to K(\bZ; 4)$, on $\pi_4$ is given by
multiplication by $b$ (up to sign; both homotopy groups are infinite cyclic). In particular, $c_2:BU \to K(\bZ;4)$ is a $\pi_4$-isomorphism.
\end{prop}

\textbf{Proof:} The quaternionic Hopf bundle on $\bS^4$ gives a map $\bS^4 \to B SU(2) \to BU$, which is a generator of $\pi_4(BU)$. The value
of the characteristic class $ac_{1}^{2} + b c_2$ on this bundle is clearly $b$.\qed

Thus we have constructed a $\pi_4$-isomorphism $\loopinf \MT \to K(\bZ;4)$. It should be remarked that for $n \geq 3$, there is no
$\pi_{2n}$-isomorphism $BU \to K(\bZ;2n)$. This is a consequence of Bott´s divisibility theorem.\\
Note that the composition $B \gaminf \to \loopinf \MT \to \loopinf \cpinf_+ \to \bZ \times BU$ is a virtual vector bundle. It is worth to
describe this vector bundle in more surface-theoretic terms. This is done in \cite{Eb}, p.69. We only state the result. Endow the universal
surface bundle $\pi:B \Gamma_{g}^{1} \to B \Gamma_g$ with a complex structure as above. Then we consider the virtual vector bundle whose fiber
over $x \in B \Gamma_g$ is given by $\overline{H^0(C_x, \omega_{C_x})}-\overline{H^0(C_x, \omega_{C_x}^{2})}$, where $C_x := \pi^{-1}(x)$,
$\omega_{C_x}$ is the canonical invertible sheaf on the Riemann surface $C_x$. This bundle is classified by the map above.

The homotopy-theoretic significance of the $\pi_4$-isomorphism above is explained by the following lemma.

\begin{lem}
Let $X$ be a simple space. Then there exists a $\pi_n$-isomorphism $X \to K(\pi_n(X);n)$ if and only if the Postnikov invariant $k^n$ is
trivial.
\end{lem}

\textbf{Proof:} The homopty class $a_n:X \to \tau_{\leq n} X$ has the following universal property. If $Y$ is an $n$-coconnected space and $f:X
\to Y$ a homotopy class, then there exists a unique homotopy class $g:  \tau_{\leq n} X \to Y$ such that $g \circ a_n = f$. Thus the existence
of a $\pi_n$-isomorphism $f:X \to K(\pi_n(X);n)$ implies the existence of a $\pi_n$-isomorphism $g:\tau_{\leq n}X \to K(\pi_n(X);n)$. In other
words, the fibration $p_n:\tau_{\leq n}X \to \tau_{\leq n-1}X$ is trivial. The classifying map $k^n$ of $p_n$ is then also trivial. Conversely,
if $k^n=0$, then $p_n$ is trivial and there is a $\pi_n$-isomorphism $g:\tau_{\leq n}X \to K(\pi_n(X);n)$ and the composition $g \circ a_n$ is
the desired $\pi_n$-isomorphism.\qed

In our case, it follows that $\tau_{\leq 4} B \gaminf^+ \simeq \tau_{\leq 3} B \gaminf^+ \times K(\bZ;4)$. To determine the homotopy type
$\tau_{\leq 4} B \gaminf^+$, we need to determine the Postnikov invariant $k^3: \tau_{\leq 2}B \gaminf^+ \to K(\bZ/24;4)$. This will occupy the rest of the paper.

\section{The icosahedral group and an interesting action of it on a surface}\label{icosahedral}

Consider a regular icosahedron $\cI$ in Euclidean $3$-space, centered at $0$ and such that all verticaes lie on $\bS^2$. It has $20$ faces (which are triangles), $12$ vertices (at every
vertex, exactly $5$ edges and $5$ faces meet) and $30$ edges. Let $G \subset SO(3)$ be the symmetry group of the icosahedron. We choose a
vertex, an edge and a face of the icosahedron and denote the isotropy subgroups by $G_2$, $G_3$ and $G_5$, respectively. These groups are cyclic
of order $2,3,5$, respectively and we choose generators $y_2, y_3 , y_5$ of these groups. The group $G$ acts transitively on the vertices as well
as on the edges and faces of the icosahedron and thus it has order $60$. It is well-known that $G$ is isomorphic to the groups $\bP \Sl_2(\bF_5)$ and $A_5$ (the alternating group). In particular, $G$ is perfect. Using the action of $G$ on the icosahedron, it is easy to see
that $\{y_2,y_3,y_5\}$ generates $G$. \\
Let $\hat{G}$ the preimage of $G$ under the $2$-fold covering $\bS^3 \to SO(3)$. The group
$\hat{G}$ is also perfect and its center is the same as the kernel of the map $\phi:\hat{G} \to G$ and contains exactly one nontrivial element
$h$. It can be shown that the extension $\bZ/2 \to \hat{G} \to G$ is the universal central extension of $G$.
The group $\hat{G}$ is the \emph{binary icosahedral group}. By the way, $\hat{G} \cong \Sl_2(\bF_5)$.\\
We choose preimages $x_i$ of the generators $y_i$. These set $\{x_2,x_3,x_5\}$ generates $\hat{G}$. The
quotient $ \hat{G} \backslash \bS^3 $ is the famous Poincar\'e homology $3$-sphere.\\
Our analysis of $G$ and $\hat{G}$ begins with the description of the Sylow-subgroups. Clearly $G_{(3)} \cong \langle y_3 \rangle \cong \bZ/3$ is a $3$-Sylow subgroup and $G_{(5)} \cong \langle y_5 \rangle \cong \bZ/5$ is a $5$-Sylow subgroup of $G$. In $\hat{G}$, $\langle x_{3}^{2} \rangle $ and $\langle x_{5}^{2} \rangle $ are Sylow subgroups.\\
The subgroup $G_{(2)}$ has order $4$ and is thus abelian. If it were cyclic, there would be $z \in G$ with $z^2=y_2$, which is obviously impossible. Thus $G_{(2)} \cong \bZ/2 \times \bZ/2=: V_4$, the generators being $y_2$ and another element which fixes a perpendicular edge of $\cI$. The preimage $\hat{G}_{(2)} = \phi^{-1}(G_{(2)}) \subset \hat{G}$ is a Sylow-$2$-subgroup and it contains a unique involution, because $\bS^3$ has
a unique involution. It follows that $\hat{G}_{(2)}$ is conjugate to $Q_8 = \{\pm 1, \pm i, \pm j , \pm j\} \subset \bS^3 \subset \bH $.
Write $V_4 = \{ 0, z_1, z_2, z_1+ z_2\}$. The homomorphism $Q_8 \to V_4$ which is defined by defined by $-1 \mapsto 0$, $i \mapsto z_1$ and $j \mapsto z_2$, is an abelianization.\\

\begin{prop}
There is an isomorphism of rings $H^*(BQ_8;\bZ) \cong \bZ [a,b,u]/I$, where $I$ is the ideal generated by $8u$, $2a$, $2b$, $ab$, $a^2$ and $b^2$.
The class $u$ has degree $4$ and it is the second Chern class of the representation $Q_8 \to \SU(2)$, while $a$ and $b$ are the first Chern classes
of the representations $Q_8\to V_4 \to \bZ /2 \subset \bS^1$ (use two different projections $V_4 \to \bZ/2$).
The cohomology ring $H^* (B \hat{G};\bZ)$ is isomorphic to $ \bZ [u] / (120 u)$, where $u \in H^4 (B \hat{G};\bZ)$ is the second Chern class of the representation $\hat{G} \to \bS^3 = \SU(2)$.\\
\end{prop}

\textbf{Proof:} Because both groups act freely on $\bS^3$, their cohomology is periodic of period $4$ and the multiplication with $u$ is a
periodicity isomorphism in both cases (use the Gysin sequence).
By Poincar\'e duality, the closed oriented manifold $\bS^3 / Q_8$ has $H^2 (\bS^3 / Q_8) \cong V_4$, $H^3(\bS^3 / Q_8) \cong \bZ$. All other cohomology in positive degrees is trivial. Now consider the commutative diagram of fibrations
\xymatrix{
\bS^3 \ar[r] \ar[d]   &  \bS^3 / Q_8 \ar[d] \\
E \bS^3 \ar[r] \ar[d] &  BQ_8  \ar[d] \\
B \bS^3 \ar@{=}[r]    &  B \bS^3 . \\}

The cohomology ring of $B\bS^3$ is well known: it is the polynomial ring, generated by the second Chern class of $B \bS^3$. The upper horizontal map has degree $8$, and these facts together with the Leray-Serre spectral sequence quickly give the first part of the proposition. \\
The calculation for $\hat{G}$ is analogous and uses that $\bS^3 / \hat{G}$ is an integral homology sphere.\qed. 

It follows from this computation that the restriction maps $H^4(B \hat{G})
\to H^4(B Q_8) \cong \bZ/8$, $H^4(B \hat{G})
\to H^4(B \bZ/3) \cong \bZ/3$ and $H^4(B \hat{G})
\to H^4(B \bZ/5) \cong \bZ/5$ are all surjective. Thus their sum $H^4 (B \hat{G}) \to H^4(BQ_8) \oplus H^4(B\bZ/3)
\oplus H^4(B \bZ/5) \cong \bZ/8 \oplus \bZ /
3 \oplus \bZ/5 \cong \bZ/120$ is an isomorphism.\\
The last fact about these groups we shall need concerns the representations of $Q_8$. The abelian group $V_4$ has four nonisomorphic
onedimensional irreducible representations and their pullbacks to $Q_8$ are still nonisomorphic (they have different characters). There are
precisely $5$ conjugacy classes in $Q_8$, and a well-known formula (see \cite{Serre}) tells us that there must be another, necessarily
$2$-dimensional irreducible representation of $Q_8$. Clearly, the representation $U:=(Q_8 \to SU(2))$ is such a representation. It is the only
one on which the central element in $Q_8$ does not act as the identity, but by $-1$. It follows: If $W$ is a $Q_8$-representation, on which the
central element acts by $-1$, then $V$ is a direct sum of copies of $U$.

\subsection{Surfaces with an action of the icosahedral group}\label{kleinsaction}

Now we construct certain actions of the binary icosahedral group on Riemann surfaces which in the end will give interesting elements in $\pi_3(B
\Gamma_{g}^{+})$. The construction is based on an easy lemma. Let $\cO(k)$ be
the $k$th tensor power of the Hopf bundle on $\bC \bP^1$. Recall that the global holomorphic sections of $\cO(k)$ can be identified with the
vector space of homogeneous polynomials of degree $k$ on $\bC^2$. Further, $\cO(k)$ is an $\Sl_2(\bC)$-equivariant bundle over the
$\Sl_2(\bC)$-space $\bC \bP^1$. If $k$ is odd, then the central element $- 1 \in \Sl_2(\bC)$ acts as $- 1$ on $\cO(k)$ (i.e. on any fiber of the bundle). If
$k$ is even, then $-1$ acts trivially on $\cO(k)$, i.e. the action descends to an action of $\bP \Sl_2(\bC)$.

\begin{lem}
Let $G \subset \bP \Sl_2(\bC)$ be a finite subgroup and let $\hat{G} \subset \Sl_2(\bC)$ be the extension of $G$ by $\bZ/2$. Let $m \in \bN$ be positive. Let $s$ be a $G$-invariant holomorphic section of $\cO(2m)$ having only simple zeroes. Then there exists a connected Riemann surface $F$ with a $\hat{G}$-action and a two-sheeted $\hat{G}$-equivariant branched covering $f:F 
\to \bC \bP^1$, which branched precisely over the zeroes of $s$.\\
If $m$ is odd, then the central element $h \in \hat{G}$ is the hyperelliptic involution on $F$, if $m$ is even, then $h$ acts trivially on $F$.\\
\end{lem}

\textbf{Proof:} Let $S \subset \cO(2m)$ be the graph of the section $s$. It is a surface of genus $0$ and it is stable under the $G$-action on
$\cO(2m)$. Let $q: \cO(m) \to \cO(2m)$ be the squaring map, let $F:=q^{-1}(S)$ and let $f:= q|_{F}$. Clearly, $F$ has a $\hat{G}$-action and $f$
is equivariant. Also, a generic point of $S$ has exactly $2$ preimages, being permuted by the involution $-1$ on the fibers on $\cO(m)$.
We need to show that $F$ is a smooth connected Riemann surface. The smoothness of $F$ is equivalent to the condition that the zeroes are simple (the locus $\{(x,y) \in \bC^2| y^2 = x^m\}$ is smooth if and only if $m=1$), and the connectivity is clear since $s$ has precisely $2m >0$ zeroes.\\
We have seen that the antipodal map of $\cO(m)$ induces the hyperelliptic involution.\qed\\

From now on, let $G$ again be the icosahedral group.

\begin{expl}\label{1}
Let $z_1, \ldots , z_{30} \in \bS^2$ be the midpoints of the edges of the icosahedron and consider them as points on $\bC \bP^1$ after the choice of a conformal map $\bS^2 \cong \bC \bP^1$. The
precise value of the points does not play a significant role in this discussion.\\
Now we take holomorphic sections $s_i$, $i=1, \ldots ,30$ of the Hopf bundle $\cO(1)$ (alias square root of the tangent bundle of $\bC \bP^1$)
having a simple zero at $z_i$ and being nonzero elsewhere. Such sections exist and are unique up to multiplication with a
complex constant. Set $s := s_1 \otimes \ldots \otimes s_{30} \in H^0(\bC \bP^1, \cO(30))$.\\
For $g \in \hat{G}$, there exists a $c(g) \in \bC^{\times}$ with $g s = c(g) s$, because $gs$ has the same zeroes as $s$. The map $c:g \mapsto
c(g)$ is a homomorphism $\hat{G} \to \bC^{\times}$. Since $\hat{G}$ has no Abelian quotient, $c$ is constant. Thus, $s$ is
an invariant section.\\
If we apply the construction of the lemma to $s$, we obtain a surface of genus $14$ (by the Riemann-Hurwitz formula) with a $\hat{G}$-action.
\end{expl}

\begin{prop}\label{numbfix}
The number of fixed points of the elements $x_2,x_3,x_5 \in \hat{G}$ in Example \ref{1} is $2, 0,0$, respectively. All powers $x_{3}^{2r}$,
$r \not \equiv 0 \pmod 3$ and $x_{5}^{2r}$, $r \not \equiv 0 \pmod 5$ have precisely $4$ fixed points.\\
\end{prop}

\textbf{Proof:}\\
Since $h$ is the hyperelliptic involution, it has precisely $2g+2 =30$ fixed points, namely the branch points which lie over the midpoints
of the edges of the icosahed\-ron.\\
There are two fixed points of $y_2$ on $\bP^1$, and they are two opposite branch points (i.e. midpoints of edges). Thus, $x_1$ has exactly
2 fixed points on $F$.\\
$y_3$ fixes precisely two opposite midpoints of faces. Hence all fixed points of $x_3$ lie over these two points. But if one of them would be a
fixed point of $x_3$, then this must also be a fixed point of $x_{3}^{3} =h$. This is impossible, since $h$ has already the $30$ fixed points
mentioned above and since no nontrivial automorphism of a Riemann surface of genus $g \leq 2$ can have more than $2g+2$ fixed points
(\cite{FaKr}, p.257). Thus $x_2$ is fixed-point-free.\\
The same argument shows that $x_3$ is also fixed-point free.\\
The second sentence follows from the same argument. \qed

Of course, one can apply a similar construction to the midpoints of the faces or to the vertices of the icosahedron. The results are surfaces of
lower genus and elements in $\pi_3(B \Gamma_g)^+$ of lower order.

\subsection{The proof of Proposition \ref{order24}}

Consider the surface $F$ of genus $14$ with the action of $\hat{G}$ constructed in the last section. Let $\rho: \hat{G} \to \Gamma_{14}$ be defined by this action. The action is holomorphic by construction and therefore, there is an induced linear action $\hat{G} \actson H^1(F; \cO)$. The cohomology class $B \rho^{*} \zeta_2$ of \ref{order24} is the same as the second Chern class of this linear $\hat{G}$-representation, which will be briefly denoted by $c$.
By the Dolbeault theorem, the following isomorphism of complex
$\hat{G}$-representations holds:

\begin{equation}\label{decomposi}
H^1(F;\cO) \otimes_{\bR} \bC \cong H^1(F; \bC).
\end{equation}

According to the structure of $H^4(B\hat{G})$ described above, we need to show: The restriction of $c$ to $H^4(B \bZ/3)$ has
order $3$, the restriction to $H^4(\bZ/5)$ is trivial and the restriction to $H^4( B Q_8)$ has order $8$.\\
Denote by $L_r$ be the irreducible representation $ 1 \mapsto \exp (\frac{2 \pi i r}{p})$ of $\bZ/p$ on $ \bC$. The Lefschetz fixed point
formula, applied to the result of Proposition \ref{numbfix}, allows us to determine the decomposition of the representation of $\bZ/3 \cong \hat
{G}_{(3)} = \langle x_{3}^{2} \rangle$ on the $28$-dimensional space $H^1(F; \bC)$. The result is $H^1(F; \bC) \cong 8 \bC \oplus 10 L_1 \oplus
10 L_2$. By \ref{decomposi}, it follows that $H^1(F; \cO) \cong 4 \bC \oplus 5 L_1 \oplus 5 L_2$ as $\bZ/3$-modules. The Chern polynomial of
$L_n$ is $1+ nv$; $ v \in H^2(\bZ/3)$ an appropriate generator. It follows that $c|_{B \bZ/3} = v^2$.\\
The result for the subgroup $\bZ/5$ follows by an analogous argument. The decomposition of $H^1(F; \cO)$ as a $\bZ/5$-module is
$2\bC \oplus 3(L_1 \oplus L_2 \oplus L_3 \oplus L_4)$.\\
We have seen that the central element $h \in Q_8 \subset \hat{G}$ acts as a hyperelliptic involution. Thus it acts on $H^1(F; \cO)$ by $-1$. As
we have seen, this implies that $H^1(F, \cO)$, as a $Q_8$-module, decomposes into seven copies of the twodimensional representation $U$. The
first Chern class of $U$ is zero, while the second is a generator of $H^4(B Q_8)$. Thus $c|_{B Q_8}$ is seven times a
generator and thus again a generator. This finishes the proof of Proposition \ref{order24}

\section{Conclusion}\label{end}

\begin{lem}
The order of $k^3$ in $H^4(K(\bZ;2); \bZ/24)$ is the quotient $\frac{24}{\sharp \Tors H^4 (B \gaminf; \bZ)}$.
\end{lem}

\textbf{Proof:} The triviality of $k^4$ shows that $H^4 (B \gaminf; \bZ)= H^4 (\tau_{\leq 4}B \gaminf; \bZ) = H^4 (\tau_{\leq 3}B \gaminf \times
K(\bZ;4); \bZ)= H^4 (\tau_{\leq 3}B \gaminf) \oplus \bZ$.

Now consider the pullback-diagram of homotopy-fibrations

\begin{equation}
\xymatrix{K(\bZ/24;3) \ar[r] \ar[d]             & K(\bZ/24;3) \ar[d]\\
\tau_{\leq 3} B \gaminf^+ \ar[r] \ar[d] &   \ast \ar[d]\\
K(\bZ;2)   \ar[r]^{k^3}                &  K(\bZ/24;)\\}
\end{equation}

and the associated Leray-Serre spectral sequences in homology. By the universal coefficient theorem, $\sharp \Tors H^4 (B \gaminf; \bZ) = \sharp
H_3 (B \gaminf^+; \bZ)$. Also, the order of $k^3$ in $H^4(K(\bZ;2); \bZ/24)$ is equal to the order of the image of $k^{3}_{*}: H_4 (K(\bZ;2);
\bZ) \cong \bZ \to H_4(K(\bZ/24;4);\bZ) \cong \bZ/24$. One deduces a commutative square

$$\xymatrix{
H_3(K(\bZ/24;3); \bZ) \ar@{=}[r] & H_3(K(\bZ/24;3); \bZ)\\
H_4(K(\bZ;2) ; \bZ) \ar[r]^{k^3} \ar[u]^{d_3} & H_4 (K(\bZ/24; 4); \bZ) \ar[u]^{d_3},
 }$$

where the $d_3$ on the right-hand-side is an isomorphism. This shows that $\sharp (H_3 (\tau_{\leq 3} B \gaminf^+; \bZ)) = \sharp (\bZ/24 / \im k^3) = \frac{24}{\ord  k^3}$. \qed

To conclude our argument, we consider the fibration

\begin{equation}
\xymatrix{ \tau_{\geq 3} B \gaminf^+ \ar[r]^-{\iota_2} & B \gaminf^+ \ar[r]^-{a_2} & \tau_{\leq 2} B \gaminf^+ = K(\bZ;2).  }
\end{equation}

The Leray-Serre spectral sequence yields a short exact sequence

\begin{equation}\label{lerrserr}
\xymatrix{ 0 \ar[r] & H^4(K(\bZ;2)) \ar[r]^-{a_{2}^{*}} & H^4(B \gaminf^+) \ar[r]^-{\iota_{2}^{*}} & H^4 (\tau_{\geq 3} B \gaminf^+) \ar[r] & 0.
}
\end{equation}

The torsion subgroup $\Tors H^4 (\tau_{\geq 3} B \gaminf^+)$ is isomorphic to $\bZ/24$ because of Theorem \ref{homotopy}. Because $\gamma_2 =
\gamma_{1}^{2} - 2 \zeta_2$ is a torsion element and because $\iota_{2}^{*} \gamma_1=0$, we see that $\iota_{2}^{*} \zeta_2 \in \Tors H^4
(\tau_{\geq 3} B \gaminf^+)$. By obstruction theory, there is a unique (up to homotopy) lift $\tilde{B \rho}:B \hat{G} \to \tau_{\geq 3} B
\gaminf^+$ of $B \rho:B \hat{G} \to B \gaminf^+$ (recall that $H^*(B \hat{G})=0$ for $0 < * < 4$). But by Proposition \ref{order24}, $\tilde{B
\rho}^* \iota_{2}^{*} \zeta_2= B \rho^* \zeta_2$ has order $24$. Thus $\iota_{2}^{*} \zeta_2$ is a generator of $ \Tors H^4 (\tau_{\geq 3} B
\gaminf^+)$. The functor which assigns to an abelian group $A$ its torsion subgroup $\Tors A$ is left-exact. Thus $\Tors H^4(B \gaminf^+) \to
\Tors H^4( \tau_{\geq 3}  B \gaminf^+)$ is injective. The element $\gamma_2 \in \Tors H^4(B \gaminf^+)$ has order $12$ or $24$, again by
Proposition \ref{order24}. Finally, we show that $\iota_{2}^{*} \zeta_2$ does not lie in the image of $\iota_{2}^{*}: \Tors H^4(B \gaminf^+) \to
\Tors H^4( \tau{\geq 3}  B \gaminf^+)$. In view of the exact sequence \ref{lerrserr}, all preimages of $\iota_{2}^{*} \zeta_2$ in $H^4(B
\gaminf^+)$ are of the form $n \zeta_{1}^{2} + \zeta_2$. If $n \zeta_{1}^{2} + \zeta_2$ is torsion, then $(2n+1) \zeta_{1}^{2}$ would be torsion
(because $\zeta_{1}^{2} - 2 \zeta_2$ is torsion), which is absurd because of \ref{lerrserr}. This finishes the proof of Theorem \ref{homology}.\\
Finally, we determine the order of $\theta \in \pi_3(B \gaminf^+)$ given by $B \rho$ as above. Consider the Bockstein sequence

$$\xymatrix{
0 \ar[r] & H^3 (\tau_{\geq 3} B \gaminf^+; \bZ/24)\ar[r]^-{\delta} & H^4 (\tau_{\geq 3} B \gaminf^+; \bZ) \ar[r]^-{24} & H^4(\tau_{\geq 3} B
\gaminf^+; \bZ) }$$

There is a unique $\beta \in H^3 (\tau_{\geq 3} B \gaminf^+; \bZ/24)$ such that $\delta (\beta)= \zeta_2$. It is a generator and thus the map
$\tau_{\geq 3} B \gaminf^+ \to K(\bZ/24;3)$ associated to $\beta$ is a $\pi_3$-isomorphism. Hence the following construction gives an isomorphism
$\psi:\pi_3 (B \gaminf^+) \to \bZ/24$. Let $f: \bS^3 \to B \gaminf^+$ and denote by $\tilde{f}: \bS^3 \to \tau_{\geq 3} B \gaminf^+$ the unique
lift. Set $\psi([f]):= \langle \tilde{f}^* (\beta); [\bS^3] \rangle$ (take the Kronecker product with $\bZ/24$-coefficients).

\begin{lem}
The map $M \to B \hat{G}$ and thus its plus-construction $\bS^3 \to B \hat{G}^{+}$ induces an isomorphism $H^3(B \hat{G}; \bZ/24) \to H^3(M,
\bZ/24)$.
\end{lem}

\textbf{Proof:} We have seen that $H_3(B\hat{G}) \cong \bZ/120$ and using similar arguments as in section \ref{icosahedral}, one can see that
$H_3(M; \bZ) \cong \bZ \to H_3(B \hat{G}; \bZ)$ is surjective. By the universal coefficient theorem, it follows that the induced map in
$\bZ/24$-cohomology is an isomorphism.\qed

Consider the composition $f: \bS^3 \to B \hat{G}^{+} \to \tau_{\geq 3} B \gaminf^+$. It is an easy diagram chase using the Bockstein sequences
for $\tau_{\geq 3} B \gaminf^+$ and $B \hat{G}^{+}$ to conclude that $H^3(f; \bZ/24)$ is an isomorphism. This shows that the class $\theta \in \pi_3(B \gaminf^{+})$ is a generator, as asserted in Theorem \ref{generator}.

\end{document}